\documentclass[11pt]{amsart}
\newtheorem{theorem}{Theorem}
\newtheorem{fact}{Fact}
\newtheorem{lemma}{Lemma}
\newtheorem{ansatz}{Ansatz}
\newtheorem{proposition}{Proposition}
\newtheorem{remark}{Remark}
\newtheorem{corollary}{Corollary}
\usepackage{amssymb}
\usepackage{latexsym}
\usepackage{graphicx}

\begin{document}
%\setlength{\baselineskip}{.20 in}
%\setlength{\parskip}{.16 in}
%\setlength{\parindent}{.25 in}

%\title{On the number of zeros of 
%certain rational harmonic functions}
%\author{Dmitry Khavinson\footnote{
%\textit{2000 Mathematics Subject Classification}.
%Primary 26C15;
%Secondary 30D05, 83C99.
%The first author
%was supported by a grant from the National Science Foundation.}
%\ and Genevra Neumann}
%\date{January 14, 2004, revised March 2, 2004}
%\maketitle
\title{On the number of zeros of certain rational harmonic functions}
\author{Dmitry Khavinson}
\address{
Department of Mathematical Sciences, 
University of Arkansas,
Fayetteville, Arkan\-sas 72701}
\email{dmitry@uark.edu}
\author{Genevra Neumann}
\address{
Department of Mathematics,
Kansas State University,
Manhattan, Kansas 66506}
\email{neumann@math.ksu.edu}
\thanks{The first author was supported by a grant from the National 
Science Foundation.}
\date{January 14, 2004, revised March 2, 2004}

\begin{abstract}
Extending a result from~\cite{ks:polys},
we show that the rational harmonic function
$\overline{r(z)} - z$, 
where $r(z)$ is a rational function of degree $n > 1$,
has no more than
$5n - 5$ complex zeros.
Applications to gravitational lensing are discussed.
In particular, 
this result settles a conjecture~\cite{r:conj}
concerning the maximum number of lensed images due to an $n$-point
gravitational lens.
\end{abstract}
\maketitle

\section{Introduction}

A. Wilmshurst~\cite{wil:polys} showed that there is an upper bound
on the number of zeros of a harmonic polynomial 
$f(z) = p(z) - \overline{q(z)}$,
where $p$ and $q$ are analytic polynomials of different degree,
answering the question of
T. Sheil-Small~\cite{ss:ques}.
Let $n = \deg\,p > \deg\,q = m$.
Wilmshurst showed that $n^2$ is a sharp upper bound when $m = n - 1$
and conjectured that the upper bound is actually $m (m - 1) + 3 n - 2$.
D. Khavinson and G. \'{S}wi\c{a}tek~\cite{ks:polys} showed
that Wilmshurst's conjecture holds for the case $n > 1, m = 1$
using methods from complex dynamics.
When hearing of this result,
P. Poggi-Corradini asked whether this approach can be extended to 
the case $f(z) = p(z)/q(z) - \overline{z}$,
where $p$ and $q$ are analytic polynomials.
 
In this note,
we apply the approach from~\cite{ks:polys} 
to prove 
\begin{theorem}
\label{thm:mainresult}
Let $r(z) = p(z) / q(z)$ be a rational function where
$p$ and $q$ are relatively prime, analytic polynomials and
such that
$n = \deg\ r = \max(\deg\ p, \deg\ q) > 1$.
Then
\begin{center}
$\#\{z \in \mathbb{C}: \overline{r(z)} = z\} \le 5n - 5$
\end{center}
\end{theorem}

We note that the zeros of $\overline{r(z)} - z$ are isolated,
because each zero is also a fixed point of 
$Q(z) = \overline{r(\overline{r(z)})}$,
an analytic rational function of degree $n^2$.
This also follows from a result of
P. Davis~\cite{da:schwarz} (Chapter 14)
concerning the Schwarz functions of analytic curves.
(The Schwarz function is an analytic function $S(z)$ that gives
the equation of a curve in the form $\overline{z} = S(z)$,
\textit{cf.}~\cite{da:schwarz}.)
A rational Schwarz function implies that the curve
is a line or a circle,
so the degree must be one.

We also note that $\overline{r(z)} - z$ will not have a zero at $\infty$.
If $\infty$ were a zero,
then $r(z) = a z + b + O(1/z)$.
By the change of variable
$z = 1/w$,
we see that 
$a/w - 1/\overline{w} + b + O(w)$
has a zero at $w = 0$. 
Restricting $w$ to the real axis,
we see that $a = 1$ and $b = 0$.
We obtain a contradiction when we restrict $w$ to the imaginary axis.

We now discuss an application of our result to gravitational microlensing.
An $n$-point gravitational lens can be modeled as follows.
Suppose that we have $n$ point masses 
(such as stars).
Construct a plane 
through the center of mass of these point masses,
such that the line of sight from the observer to 
the center of mass is orthogonal to this plane.
This plane is called the lens plane
(or deflector plane).
Suppose that the lens plane is between the observer and a light source.
(We are assuming that the distance between the point masses 
is small compared to the
distance between the observer and the lens plane,
as well as the distance between the lens plane and the light source.)
The plane containing our light source which is parallel to the lens
plane is called the source plane.
Due to the deflection of light by the point masses,
multiple images of the light source are formed.
This phenomenon is known
as gravitational microlensing.
See~\cite{w:lensing} and~\cite{nb:lensing} 
for an introduction to gravitational lensing and
\cite{s:survey} for an introduction to a complex formulation of
lensing theory;
also see~\cite{plw:book}.

Gravitational microlensing can be modeled by a lens equation,
which defines a mapping from the lens plane to
the source plane.
To set up a lens equation for our $n$-point gravitational lens,
the point masses are projected onto the lens plane.
The projection of the $j$-th point mass 
has a scaled mass of $m_j$ and is located at 
a scaled coordinate of $z_j$ in the lens plane,
where $m_j$ is a positive constant and $z_j$ is a complex constant.
Suppose that we have a light source located at a scaled
coordinate of $w$ in the source plane.
Following \cite{w:micro},
this lens equation will be given by
$w = z\, + \gamma \overline{z}\, -\, 
$sign$(\sigma) \Sigma_{j=1}^n\, m_j / (\overline{z}\, -\, \overline{z_j})$,
where the normalized (external) shear $\gamma$ and 
the optical depth (or normalized surface density)
$\sigma \neq 0$
are real constants.
In this model,
if $z$ satisfies the lens equation,
then our gravitational lens will map $z$ to $w$;
hence $z$ corresponds to the position of a lensed image.
The number of lensed images
is the number of solutions of the lens equation.

We can rewrite this lens equation in terms of the rational harmonic function
$f(z)= \overline{r(z)} - z$ by letting
$r(z) = \overline{w}\,- \gamma z\, +\, 
$sign$(\sigma) \Sigma_{j=1}^n\, m_j / (z - z_j)$.
We thus see that the zeros of 
$f(z)$ are solutions of the lens equation for a 
light source at position $w$. 
H. Witt~\cite{w:micro} showed for $n > 1$ that the maximum number of 
observed images
is at most $n^2 + 1$
when $\gamma = 0$
and $(n + 1)^2$
when $\gamma \neq 0$.
A. Petters~\cite{p:morse} used Morse theory to obtain further estimates for
both cases.
S. Mao, A. Petters, and H. Witt~\cite{mpw:3n+1} 
conjectured that the maximum number
of images is linear in $n$
for the case $\gamma = 0$ and $\sigma > 0$.
S. H. Rhie~\cite{r:conj} later conjectured that for $n > 1$ such a 
gravitational lens gives at most $5 n - 5$ images.
In the $\gamma = 0$ case,
$\deg\,r = n$ and
Theorem 1 settles this conjecture.
Further,
for the case $\gamma \neq 0$,
we see that $\deg\,r = n + 1$,
so Theorem 1 gives an upper bound of $5 (n + 1) - 5 = 5 n$ lensed images.
As a result,
we have the following corollary:

\begin{corollary}
An $n$-point gravitational lens modeled by
$w = z\, + \gamma \overline{z}\, -\, 
$sign$(\sigma) \Sigma_{j=1}^n\, m_j / (\overline{z}\, -\, \overline{z_j})$,
where $n > 1$
can produce at most $5n - 5$ images
when the shear $\gamma = 0$ and at most $5n$ images when the shear
$\gamma \neq 0$.
\end{corollary}

\section{Acknowledgements}
The authors are indebted to Professor Jeffrey Rabin of the
University of California, San Diego,
for pointing out that our result settles a conjecture in gravitational
lensing,
for alerting us to the papers of S. H. Rhie 
(\cite{r:conj} and \cite{r:sharp}),
and for noting that~\cite{r:sharp} addresses the question of sharpness.
We wish to thank Professor Donald Sarason of the University of California,
Berkeley for helpful discussions concerning Rhie's sharpness proof.
The second author is grateful to Professor Larry Weaver of 
Kansas State University for 
helpful discussions concerning 
gravitational lensing.
Both authors wish to thank 
Professor Shude Mao of the Jodrell Bank Observatory at the 
University of Manchester
and Professor Arlie Petters of Duke University for valuable comments
on the second draft of this paper.

\section{Preliminaries}

We first recall some terminology and a fact from complex dynamics.
Let $r(z) = p(z) / q(z)$ be a rational function,
where $p$ and $q$ are relatively prime, analytic polynomials.
The degree of $r$, 
denoted by $\deg\,r$, 
is given by $\max(\deg\,p,\deg\,q)$.
$z$ is called a \textit{critical point of} $r$ if $r'(z) = 0$.
We will be interested in counting the fixed points of a rational
function.
A fixed point $z_0 \in \mathbb{C_\infty}$ is said to be
\textit{attractive}, \textit{repelling} or \textit{neutral} if
$|r'(z_0)| < 1$, $|r'(z_0)| > 1$, or $|r'(z_0)| = 1$
respectively.
A neutral fixed point where the derivative is a root of unity is said to
be \textit{rationally neutral}.
A fixed point $z_0$ is said to \textit{attract} some point 
$w \in \mathbb{C_\infty}$
if the iterates of $r$ at $w$ converge to $z_0$.
We will be using the following result,
whose proof can be found in~\cite{cg:dynamics}, 
Chapter III, Theorems 2.2 and 2.3:

\begin{fact}
\label{fact:fatou}
Let $r$ be a rational function with $\deg\,r > 1$.
If $z_0$ is an attracting or rationally neutral fixed point,
then $z_0$ attracts some critical point of $r$.
\end{fact}

We will also use a version of the argument principle for harmonic functions.
Let $f$ be harmonic in an open set. 
We say that $z$ is in the \textit{critical set} of $f$ if
the Jacobian of $f$ vanishes at $z$.
In particular,
in a sufficiently small disk containing $z$,
there exist analytic functions $h$ and $g$ such that $f = h + \overline{g}$
in this disk.
Then $z$ is in the critical set of $f$ if $|h'(z)|^2 - |g'(z)|^2 = 0$.
We say that $z$ is a \textit{singular zero} of $f$ if $z$ is in the
critical set of $f$ and $f(z) = 0$.
$f$ is said to be \textit{sense-preserving} at $z$ if
$|h'(z)| > |g'(z)|$
and is said to be \textit{sense-reversing} at $z$ if
$|h'(z)| < |g'(z)|$.
The \textit{order} of a non-singular zero is positive if $f$ 
is sense-preserving at the zero and negative if $f$ is sense-reversing.

Let $f$ be harmonic in a punctured neighborhood of $z_0$.
We will refer to $z_0$ as a \textit{pole} of $f$
if $f(z) \rightarrow \infty$ as $z \rightarrow z_0$.
Let $C$ be an oriented closed curve that contains neither zeros
nor poles of $f$.
The notation $\Delta_C \arg\,f(z)$ denotes the increment in the
argument of $f(z)$ along $C$.
Following~\cite{st:part},
the \textit{order} of a pole of $f$ is given by 
$-\frac{1}{2 \pi} \Delta_C \arg\,f(z)$,
where $C$ is a sufficiently small circle around the pole.
We note that if $f$ is sense-reversing on a sufficiently small circle
around the pole,
then the order of the pole will be negative.
We will use the following version of the argument principle 
which is taken from~\cite{st:part}:
\begin{fact}
\label{fact:argprinc}
Let $f$ be harmonic, 
except for a finite number of poles,
in a simply connected domain $D$ in the complex plane.
Let $C$ be a Jordan curve contained in $D$ not passing through a pole
or a zero,
and let $\Omega$ be the open bounded region created by $C$.
Suppose that $f$ has no singular zeros in $D$ and let $N$ be the sum
of the orders of the zeros of $f$ in $\Omega$.
Let $M$ be the sum of the orders of the poles of $f$ in $\Omega$.
Then $\Delta_C \arg\,f(z) = 2 \pi N - 2 \pi M$.
\end{fact}
\noindent
We note that a more general form of the argument principle 
can be found in~\cite{ss:polys}.

\section{Non-repelling fixed points}

\begin{proposition}
\label{prop:sensepreszeros}
Let $p$ and $q$ be relatively prime analytic polynomials.
If $r = p / q$ is a rational function of degree $n > 1$,
then the set of points for which 
$z = \overline{r(z)}$ and $|r'(z)| \le 1$ has
cardinality at most $2n - 2$.
\end{proposition}

\noindent
\textit{Proof.}
Let $n_+$ denote the number of points $z_0$ satisfying the conditions
of Proposition \ref{prop:sensepreszeros}.
Following the approach of D. Khavinson and
G. \'{S}wi\c{a}tek~\cite{ks:polys},
we consider the function $Q(z) = \overline{r(\overline{r(z)})}$,
which is a rational function of degree $n^2$.
As in~\cite{ks:polys},
we note that all fixed points of $\overline{r(z)}$ that are critical points
for $f(z) = \overline{r(z)} - z$ 
are rationally neutral fixed points for $Q(z)$,
so Fact \ref{fact:fatou} will apply.
Their outline (Lemmas 1 - 3 in~\cite{ks:polys} together 
with Fact \ref{fact:fatou}) 
carry over word for word to the rational case if $\mathbb{C}$ is
replaced by the Riemann sphere
$\mathbb{C}_\infty$.
In particular,
each point $z_0$ which satisfies the conditions of
Proposition \ref{prop:sensepreszeros} attracts at least $n + 1$
critical points of $Q$.

Since $\deg\ Q = n^2$,
by the Riemann-Hurwitz relation 
(see~\cite{fo:surfaces}, 
Section 17.14),
$Q$ has $2 n^2 - 2$ critical points 
(counted with multiplicities)
in $\mathbb{C}_\infty$.
$n_+$ will be largest when all of the critical points of $Q$
are attracted to points $z_0$ satisfying the conditions
of Proposition \ref{prop:sensepreszeros}.
Since different fixed points attract
disjoint sets of critical points,
we see that
$2 n^2 - 2 \ge (n + 1)\,n_+$ and the claim follows.
$\Box$

\section{Proof of Theorem \ref{thm:mainresult}}

As in~\cite{ks:polys},
we say that the rational function $r(z)$ is \textit{regular}
if no critical point of $f(z) = \overline{r(z)} - z$
is a zero of $f$.
This will allow us to apply Fact \ref{fact:argprinc} to count the zeros
of $f$.

\begin{ansatz}
\label{lem:senserevzeros}
If $r$ is regular of degree $n > 1$,
then $f(z) = \overline{r(z)} - z$ has at most $5 n - 5$ distinct
finite zeros .
\end{ansatz}

\noindent
\textit{Proof.}
Let $r = p / q$,
where $p$ and $q$ are relatively prime, analytic polynomials.
Let $n_p = \deg\ p$ and $n_q = \deg\ q$.
Then $n = \max\{n_p, n_q\}$.
Let $n_-$ denote the number of sense-reversing zeros of $f$
and let $n_+$ be defined as in the proof of Proposition 1.
Since $r$ is regular,
no zero of $f$ lies on the critical set of $f$;
this will allow us to apply Fact \ref{fact:argprinc}.
In particular,
$n_+$ counts the number of sense-preserving zeros of $f$.
Then $f$ has $n_+ + n_-$ zeros in the finite complex
plane,
counting multiplicity.
Note that $f$ is sense-reversing at each of the $n_q$
poles of $r$;
hence,
the order of each pole is negative.

Consider the increment in the argument of $f$ in a region
bounded by a circle of large enough radius so that all
of the finite zeros and poles of $f$ are enclosed.
We will apply Fact \ref{fact:argprinc} and use our bounds on $n_+$ from 
Proposition \ref{prop:sensepreszeros}.
Note that we are counting zeros with multiplicity.
We consider two cases.

\begin{enumerate}
\item 
Assume that $n = n_q \ge n_p \ge 0$.
Then the critical set of $f$ is bounded.
For $z$ large,
$r(z)$ is at most $O(1)$.
Thus $f$ is sense-preserving on our circle
with an argument change of $2 \pi$.
We also note that $\infty$ is not a zero of $f$.
By Fact \ref{fact:argprinc},
$1 = (n_+ - n_-) - (-n) \le 2 n - 2 - n_- + n$.
Hence $n_- \le 3 n - 3$,
so $f$ has at most
$5 n - 5$ zeros in $\mathbb{C}_\infty$
and the claim follows.

\item
Assume that $n = n_p \ge n_q + 1$.
Since $p$ and $q$ are relatively prime,
at least one of the two polynomials has a non-zero constant
term.

We first suppose that $p$ has a non-zero constant term.
Then $f(0) \neq 0$.
Let $z = 1/w$ and consider 
$F(w) = 1/\overline{r(1/w)} - w = \overline{G(w) / H(w)}\ -\ w$.
Then $F$ satisfies the conditions of the previous case,
replacing $n_q$ by $\deg\ H$ and $n_p$ by $\deg\ G$.
Thus,
$F$ has at most $5 n - 5$ zeros.
Since $f(0) \neq 0$,
$f$ must also have at most $5 n - 5$ zeros 
in the finite
complex plane.

We now suppose that $p$ does not have a non-zero constant term.
Since $p$ and $q$ are relatively prime, 
the lowest order term in $q$ must be a non-zero constant.
Consider $f_c(z) = \overline{(p(z) + c) / q(z)} - z
= \overline{r_c(1/w)} - 1/w$.
For all $c$ sufficiently small, 
we have that
$p(z) + c$ and $q$ are relatively prime.
In that case, 
$f$ and $f_c$ have the same poles
and $f_c$ approaches $f$ uniformly as $c \rightarrow 0$
on compact subsets of $\mathbb{C}$
which do not contain any of the poles of $f$.
We then substitute $z = 1/w$ to form 
$F_c(w) = 1/\overline{r_c(1/w)} - w$
$= \overline{G_c(w) / H_c(w)} - w$.
By construction,
$\deg\ G_c = \deg\ H_c = n$,
so $F_c$ has at most $5 n - 5$ zeros in the finite plane.
Hence, 
this also holds for $f_c$ since $f_c(0) \neq 0$.

Now suppose that $f$ has more than $5 n - 5$ zeros in $\mathbb{C}$.
By constructing a sufficiently small circle around
each zero of $f$,
we can guarantee that $f$ is harmonic in each of the
resulting closed disks. Moreover,
we can make our circles sufficiently small so that 
each closed disk contains no other zeros of $f$
and contains no critical points of $f$
(recall that $r$ is regular).
Assume that $f_c$ does not
vanish in such a disk surrounding a zero, 
so $|f_c| > a > 0$. 
Then,
for $c$ sufficiently small,
$\arg\,f = \arg\,(f_c\, (1 + (f - f_c) / f_c))$ 
has zero increment around the circle bounding this disk,
by Fact \ref{fact:argprinc}; 
this contradicts $f$ having a zero in this disk.
We see that $f_c$ must have the same number of zeros
as $f$ in that disk (counting multiplicity).
Thus,
$f_c$ must have more than $5 n - 5$ zeros in the finite plane,
a contradiction.
$\Box$
\end{enumerate}

It remains to show that it is enough to consider
regular rational functions.

\begin{lemma}
\label{lem:regdense}
If $r(z)$ is a rational function of degree greater than $1$,
then the set of complex numbers $c$ for which $r(z) - c$
is regular is open and dense in $\mathbb{C}$.
\end{lemma}

\noindent
\textit{Proof.}
As in~\cite{ks:polys},
it is enough to show that the image of the critical
set of $f(z) = \overline{r(z)} - z$ is nowhere
dense in $\mathbb{C}$.
In contrast to the polynomial case,
it is possible for the critical set to be unbounded.
This difficulty however is easily resolved by restricting $f$ to an
increasing family of concentric balls whose radii run to $\infty$.
$\Box$

Lemma \ref{lem:regdense} shows that the set of regular rational functions
is dense in the topology of uniform convergence in the spherical metric.
As at the end of the proof of the Ansatz,
it is easily seen that every zero of the function $f(z) = \overline{r(z)} - z$
must be a limit point for zeros of $\overline{r(z)} - z - c$ where
$c \rightarrow 0$.
Thus,
Lemma \ref{lem:regdense} and the Ansatz prove the theorem.

\section{Final remarks}
\begin{figure}
\begin{center}
\includegraphics[scale=.95]{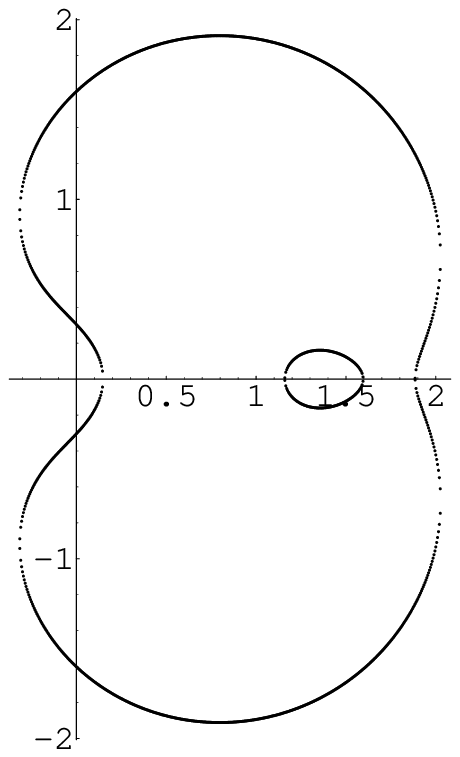}
\end{center}
\caption{Critical set for 
$f(z) = \frac{\overline{z}^2 + \overline{z} - 
\frac{1}{2}}{\overline{z}^2 - \frac{3}{2} \overline{z} + 1} - {z}$}
\label{fig:crit}
\end{figure}

\begin{figure}
\begin{center}
\includegraphics[scale=.95]{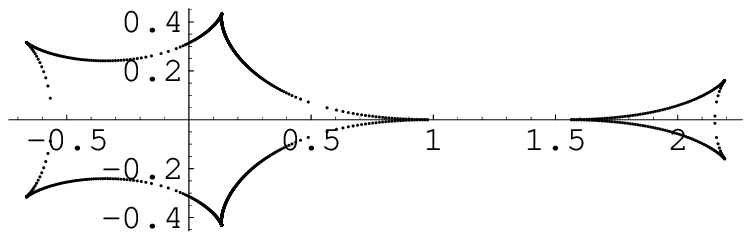}
\end{center}
\caption{Image of critical set for 
$f(z) = \frac{\overline{z}^2 + \overline{z} - 
\frac{1}{2}}{\overline{z}^2 - \frac{3}{2} \overline{z} + 1} - {z}$}
\label{fig:imag}
\end{figure}

Our proof does not indicate whether the 
$5 n - 5$ bound is always sharp.
We have, however, found an example where this bound is attained
for the case $n = 2$;
namely,
$f(z) = (\overline{z}^2 + \overline{z} - \frac{1}{2}) / 
(\overline{z}^2 - \frac{3}{2} \overline{z} + 1) - z$.
This function has $5 n - 5 = 5$ distinct finite zeros.
This function has two poles:
$z = \frac{1}{4} (3 - i\,\sqrt{7})$ and
$z = \frac{1}{4} (3 + i\,\sqrt{7})$;
$f$ is sense-reversing at these poles.
There are three sense-reversing zeros:
$z = \frac{1}{2}$, 
$\frac{1}{2}(1 + i\,\sqrt{11})$,
and
$\frac{1}{2}(1 - i\,\sqrt{11})$.
As expected from Fact \ref{fact:argprinc}
for an overall index change of $+1$ on large circles,
there are two sense-preserving zeros; 
namely,
$z = 1 - \sqrt{2}$ and $1 + \sqrt{2}$, 
Figure \ref{fig:crit} is a Mathematica plot of the critical set 
of this function and
Figure \ref{fig:imag} shows its image.
We note that $f$ is sense-preserving in the unbounded component 
in Figure \ref{fig:crit},
sense-reversing in the larger of the two bounded components,
and sense-preserving in the smaller of the two bounded components.
We also note that $f$ cannot be rewritten to model 
a 2-point gravitational lens.

L. Geyer~\cite{g:sharp} has recently shown that the $3 n - 2$ bound
on the number of zeros of $f(z) = p(z) - \overline{z}$ where
$\deg\,p = n$ is sharp for all $n > 1$.
D. Bshouty and A. Lyzzaik~\cite{bl:sharp} have recently given an
elementary proof for $n = 4, 5, 6, 8$. 
Hence,
a sharp bound on the number of zeros of $f(z) = \overline{r(z)} - z$
must be at least $3 n - 2$.

Results from gravitational lensing help with 
the question of sharpness of the $5 n - 5$ bound in Theorem 1.
We restrict our attention to the $\gamma = 0$, $\sigma > 0$ case
of the lens equation
(\textit{cf.} Introduction).
Mao, Petters, and Witt~\cite{mpw:3n+1}  have shown there are
$3 n + 1$ images when the light source is at the origin and
the $n > 2$ point masses 
(each of mass $1/n$)
are vertices of a regular
polygon centered at the origin.
Thus,
in contrast to the case where $r(z)$ is a polynomial,
by moving the poles of $r(z)$ into the finite plane,
we see that a sharp bound on the number of zeros must be at least $3 n + 1$.
Further,
by a clever perturbation argument,
Rhie~\cite{r:sharp}
has shown that the $5 n - 5$ bound is 
attained and is hence sharp for all $n > 1$.

It is known that the number of lensed images is odd in the case
of gravitational lensing by a regular
gravitational lens
(\cite{s:survey}, \cite{b:images}),
given that
the position of the source is not a critical value of the lens equation.
(A regular lens has a smooth mass distribution.
A point $w$ is a critical value of the lens equation if
the Jacobian of the lens mapping defined by the lens equation 
vanishes at any of the points in $f^{-1}(w)$.)
Given that $5n - 5$ is not always odd,
it may seem surprising that $5 n - 5$ lensed images may be possible
for an $n$-point gravitational lens with $\gamma = 0$ and $\sigma > 0$.
However,
it was shown by Petters ~\cite{p:morse} using Morse theory,
that this is not a problem.
We summarize this result in the following corollary and
present an alternate proof.

\begin{corollary}
Suppose that we have an $n$-point gravitational lens
modeled as above
with $\gamma = 0$ and a light
source which is not located at a critical value of the lens equation.
Then the number of lensed images 
must be even when $n$ is odd and
odd when $n$ is even.
\end{corollary}

\noindent
\textit{Proof.}
For this case of the lens equation,
we may apply step (1) in the proof of the Ansatz.
In particular,
the number of images $N$ will be $n_+\,+\,n_-$,
where $n_+$ denotes the number of sense-preserving solutions
and $n_-$ the number of sense-reversing solutions of the lens equation.
From step (1),
we see that $n_+ = 1 + n_- - n$,
so $N = 1 + 2 n_- - n$.
$\Box$

\bigskip
\noindent
This result has also been shown by other techniques in~\cite{r:conj} 
(footnote 2, page 2),
which extends the approach of W. Burke~\cite{b:images}
for the case of a regular gravitational lens
to that of an $n$-point gravitational lens.
Note that the above proof is similar to the approach used by 
N. Straumann~\cite{s:survey}
for the case of a regular gravitational lens.

We also note that Theorem 1 can be applied to the case of a
more general mass distribution than point masses.
For a compactly supported mass distribution 
with the projected mass density $d\psi$ in the lens plane,
the lens equation transforms into
$w = z\,+\,\gamma \overline{z}\,-\,sign(\sigma) 
\int_{\mathbb{C}} d\psi(\zeta)\,
/\,(\overline{z} - \overline{\zeta})$
(\textit{cf.}~\cite{nb:lensing} and \cite{s:survey}).
In particular,
we have the following corollary:

\begin{corollary}
Suppose that the projected mass density of 
a gravitational lens consists of $n > 1$ 
radially symmetric, continuous, compactly supported densities in the lens plane.
If the shear $\gamma = 0$,
then the number of lensed images outside the support of the masses
cannot exceed $5n - 5$.
If the shear $\gamma \neq 0$,
then the number of lensed images outside the support of the masses
cannot exceed $5n$.
\end{corollary}

\noindent
\textit{Proof.}
Consider the integral term in the lens equation for one of the mass densities.
A calculation shows that this integral
evaluated for any $z$ outside of the support of this radially symmetric
mass density
equals $b / (z - a)$,
where $a$ is the center of the mass density and $b$ is a finite constant.
Hence,
for $z$ outside of the support of the masses,
the lens equation reduces to the lens equation for $n$ point masses
(each mass is reduced to a point mass at the sphere's center).
Corollary 1 can then be applied to count the possible images,
giving the upper bounds on images outside the support of the masses
as claimed.
$\Box$

\begin{remark}
If the mass distribution is composed entirely of luminous matter,
then the only lensed images that will be visible are those that lie
outside of the support of the masses.
However,
there are mass distributions that are not luminous
(see~\cite{apod:cluster} and \cite{apod:dark}).
For such mass distributions,
it is possible to have lensed images that lie 
inside the support of the masses.
\end{remark}

A question related 
to bounding the number of zeros of $f(z) = \overline{r(z)} - z$
would be to find a bound on the number of zeros
of $f(z) = R(z) - \overline{r(z)}$,
where both $R$ and $r$ are rational functions.
It is not clear what condition would be needed for such a function 
to have a finite number of zeros,
much less a bound on the number of distinct zeros
(beyond the obvious B\'{e}zout Theorem's bound when we assume
that there are a finite number of zeros).
For example,
let $p(z)$ be a polynomial (of degree at least two) and consider
$f(z) = p(z)\,-\,1 / \overline{p(z)}$.
The zeros of this function form a lemniscate in the complex plane.
In other words,
$f(z)$ has an infinite number of zeros
and these zeros are not isolated.
Moreover,
$\lim_{z \rightarrow \infty} f(z) = \infty$.
If we compare this with the result of
Wilmshurst~\cite{wil:polys} for the case when $f$ is entire:
$f(z) \rightarrow \infty$ as $z \rightarrow \infty$
implies a finite number of zeros, 
we see that all of the zeros for an entire function
will be isolated under Wilmshurst's condition for finite valence.
Thus,
$f$ having finite poles makes the problem of discreteness of its zero set
much more subtle.

%\bigskip
%\noindent
%Department of Mathematical Sciences, 
%University of Arkansas,
%Fayetteville, Arkansas 72701\\
%\textit{E-mail address:} dmitry@uark.edu
%
%\noindent
%Department of Mathematics,
%Kansas State University,
%Manhattan, Kansas 66506\\
%\textit{E-mail address:} neumann@math.ksu.edu
\end{document}